\documentclass[oneside]{amsart}
\usepackage{amssymb}
\usepackage{geometry}
\usepackage{microtype}
\usepackage{booktabs}
\usepackage{array}
\usepackage{enumitem}
\usepackage{threeparttable}
\usepackage{placeins}
\usepackage{pgfplots}
\pgfplotsset{compat=1.18}

\usepackage[bookmarks=false]{hyperref}
\hypersetup{
colorlinks=true,
allcolors=black
}

\newcolumntype{E}{r @{ = } l}

\title{The fourth known primitive solution to \(a^5 + b^5 + c^5 + d^5 = e^5\)}
\author{Jeffrey M. Braun}
\address{New York, NY, USA}
\email{jmbraun.research@gmail.com}
\date{August 16, 2026}

\subjclass[2020]{Primary 11D41; Secondary 11Y50}
\keywords{Diophantine equations, computational number theory, fifth powers, sum of powers}

\begin{document}

\begin{abstract}
We report the fourth known primitive solution to the Diophantine equation $a^5 + b^5 + c^5 + d^5 = e^5$, extending the small set of solutions identified in 1966, 1996, and 2004.
The new solution has target value $e = 1{,}956{,}878$, more than an order of magnitude larger than those of all
previously known primitive solutions.
It was obtained by a large-scale computational search using an optimized meet-in-the-middle algorithm.
To make this search feasible under fixed memory, we introduce a partitioning strategy that trades additional preprocessing for reduced memory usage, and show that this approach remains practical for search ranges well beyond those explored here.
Finally, we compare structural properties of known fourth- and fifth-power counterexamples to Euler's sum-of-powers conjecture, observing that most known fifth-power solutions have their largest summand close to the target value.
\end{abstract}

\maketitle
\pagestyle{plain}

\section{Introduction}

Diophantine equations of the form
\begin{equation}
a_1^n + \dots + a_{k-1}^n = a_k^n
\label{eq:sum_of_powers}
\end{equation}
originate in Euler's sum-of-powers conjecture, which posited that for $n > 3$, any nontrivial solution requires $k - 1 \geq n$.
This conjecture was famously disproved in 1966 by Lander and Parkin \cite{Lander1966, Lander1967}, who identified a four-term counterexample for $n = 5$.

The understanding of such counterexamples varies significantly with $n$.
For $n = 4$, Elkies proved the existence of infinitely many primitive solutions \cite{Elkies1988}.
Computational searches have since identified many small primitive solutions \cite{Elkies1988, Bernstein2001, Gerbicz2009}, enabling further study of their algebraic and distributional properties.

By contrast, the $n = 5$ case is far less understood.
It is unknown whether infinitely many primitive solutions exist, and only three have been discovered since 1966.
This motivates the computational search undertaken here.

We investigate the $n = 5$ case
\begin{equation}
a^5 + b^5 + c^5 + d^5 = e^5
\label{eq:main}
\end{equation}
and report a fourth primitive solution. A solution is primitive if $\gcd(a,b,c,d,e)=1$.

\section{Known primitive solutions}

Table \ref{tab:known_solutions} lists the known primitive solutions to \eqref{eq:main}.

\begin{table}[htbp]
\centering
\caption{Known primitive solutions to \eqref{eq:main}.}
\label{tab:known_solutions}
\begin{threeparttable}
\renewcommand{\arraystretch}{1.2}
\begin{tabular}{E c c r}
\toprule
\multicolumn{2}{c}{Solution} & Domain & Year \\
\midrule
$27^5 + 84^5 + 110^5 + 133^5$ & $144^5$ & $\mathbb{Z}_{\geq 0}$ & 1966\\
$(-220)^5 + 5027^5 + 6237^5 + 14068^5$ & $14132^5$ & $\mathbb{Z}$ & 1996\\
$55^5 + 3183^5 + 28969^5 + 85282^5$ & $85359^5$ & $\mathbb{Z}_{\geq 0}$ & 2004\\
\bottomrule
\end{tabular}

\begin{tablenotes}[flushleft]
\item \footnotesize \textit{Sources:} \cite{Lander1966, Scher1996, Frye2004}.
\end{tablenotes}
\end{threeparttable}
\end{table}

\FloatBarrier

\section{New solution}

We report a fourth primitive solution to \eqref{eq:main} in $\mathbb{Z}$, with one negative term:
{\setlength{\abovedisplayskip}{8pt}
\setlength{\belowdisplayskip}{8pt}
\begin{equation}
  719115^5 + 1331622^5 + (-1340632)^5 + 1956213^5 = 1956878^5.
  \label{eq:new_solution}
\end{equation}
}

Equation \eqref{eq:new_solution} and its primitivity were independently verified using arbitrary-precision arithmetic.
Its target value $e = 1{,}956{,}878$ exceeds those of all known primitive solutions by more than an order of magnitude.

A notable structural feature is the near-equality of the two largest terms, discussed in Section~\ref{sec:distribution}.

\section{Search Algorithm}

\subsection{Algorithm in \(\mathbb{Z}_{\geq 0}\)}
\label{subsec:algorithm}

To search for solutions with \(e_{\min} \le e \le e_{\max}\), the algorithm uses a meet-in-the-middle approach with two distinct phases.

In the \textit{preprocessing phase}, we form the multiset
\[
S = \big( x_1^5 + x_2^5 \big)_{0 \le x_1 \le x_2 \le e_{\max}},\footnote{
Throughout the paper, $S$ denotes a context-dependent multiset of sums $x_1^5 + x_2^5$, with one entry per index pair $(x_1,x_2)$, ranging from the full multiset to filtered submultisets.
}
\]
which is then sorted in ascending order.

In the \textit{evaluation phase}, for each \(e \in [e_{\min}, e_{\max}]\), we seek \(s_1, s_2 \in S\) such that
\[
s_1 + s_2 = e^5,\quad s_1 \le s_2.
\]
We perform a two-pointer search over $S$ by initializing $s_1 = 0$ and $s_2 = \max\{ s \in S : s < e^5 \}$, and updating the pair $(s_1, s_2)$ based on the sign of
\[
D = e^5 - (s_1 + s_2).
\]
We advance $s_1$ when $D \ge 0$ and decrement $s_2$ when $D < 0$, ensuring any pair from $S$ summing to $e^5$ remains between the current pointers.
If $D = 0$, we obtain a solution which is checked for primitivity. Iterations continue until the two pointers cross.

Setting \(L = e_{\min}\) and \(N = e_{\max}\), time complexity is \(O(N^2 \log^2 N)\) for preprocessing and \(O((N - L) \cdot N^2 \log N)\) for evaluation.
To amortize preprocessing cost across the evaluation phase, intervals should be chosen with \(N - L = \Omega(N)\), giving overall time complexity \(O(N^3 \log N)\).
The algorithm requires \(O(N^2 \log N)\) space to store \(S\).

\subsection{Comparison to Alternative Algorithms}
\label{subsec:alt_algorithms}

Two alternatives to our two-pointer search are a hash-based search and a priority-queue search.

In a hash-based search, one enumerates all pairs $(a,b)$ with $a \le b$ and stores the sums $a^5 + b^5$ in a hash table.
Each triple $(c, d, e)$ with $c \le d < e$ is then tested by checking whether $e^5 - d^5 - c^5$ appears in the table.
This search has the same asymptotics as our two-pointer search.

In a priority-queue search, one retrieves $s_1$ and $s_2$ not from a sorted multiset, but from a pair of heap-backed priority queues per target $e$, in the style of \cite{Bernstein2001}.
Relative to the two-pointer search, a priority-queue search reduces space complexity by a factor of $N$ to $O(N \log N)$, at the cost of an additional $\log N$ time factor, giving $O(N^3 \log^2 N)$ time.

In benchmarking, the two-pointer search achieved a 10-fold speedup relative to the hash-based search despite evaluating only 5\% fewer candidates, consistent with the cache-locality advantages of sequential memory access.
The two-pointer search achieved a 15-fold speedup at $e_{\max} = 2^{10}$ relative to the priority-queue search, consistent with the latter's two $O(\log N)$ heap operations per iteration.
The two-pointer search is therefore strongly preferred, provided its $O(N^2 \log N)$ memory footprint can be managed for our search ranges.

\subsection{Partitioning}
\label{subsec:partitioning}

For large $e_{\max}$, the full multiset $S$ does not fit in memory. We therefore partition the computation into independent subproblems using modular residue conditions.

Let
\[
M = p_1 p_2 \cdots p_k
\]
be a product of distinct primes such that $x \mapsto x^5 \bmod p_i$ is injective for each $p_i$ (equivalently, $\gcd(5,p_i-1)=1$).

We index subproblems by unordered residue pairs $\{r_1,r_2\}$ modulo $M$, consisting of all $(s_1,s_2)$ with
\[
\{ s_1 \bmod M,\, s_2 \bmod M \} = \{ r_1, r_2 \}.
\]
For a given target value $e$, only residue pairs satisfying
\[
r_1 + r_2 \equiv e^5 \pmod{M}
\]
can contribute to solutions. Each pair defines an independent subproblem, in which we
\begin{enumerate}[label=(\roman*)]
\item generate only sums $s \in S$ with residues in $\{r_1,r_2\}$ modulo $M$, and
\item consider only target values $e$ compatible with the pair.
\end{enumerate}

Since each pair $\{s_1,s_2\}$ belongs to exactly one residue pair $\{r_1,r_2\}$, these subproblems form a partition of the search space of candidate pairs.

For example, when $M = 5$ and $e^5 \equiv 0 \pmod{5}$, the three residue pairs $\{0,0\}$, $\{1,4\}$, and $\{2,3\}$ cover all candidate pairs and can be processed as independent subproblems.
Repeating this for each residue class of $e^5 \bmod 5$ covers the full range of target values $e$.

For a given range $[e_{\min}, e_{\max}]$, this construction yields \(M(M+1)/2\) subproblems.
Off-diagonal pairs contain approximately \(2/M\) of all sums, and diagonal pairs approximately \(1/M\), reflecting the uniform residue distribution guaranteed by injectivity of \(x^5 \bmod M\).
This bounds the memory requirement per subproblem.

\subsection{Scaling Behavior under Fixed Memory}
\label{subsec:fixed_memory}

Under a fixed memory budget, partitioning imposes a growing computational cost as $e_{\max}$ increases.
This cost arises from generating and sorting each sum $x_1^5 + x_2^5$ across all $M$ residue-compatible subproblems rather than just once.
To quantify this effect, we compute the preprocessing-to-evaluation ratio under fixed memory, and show that it grows without bound as $e_{\max} \to \infty$, in contrast to the unpartitioned case.

Consider a per-subproblem memory budget fixed at $B$ bits.
With our algorithm requiring $\Theta(N^2 \log N)$ space, the modulus $M$ required to satisfy this budget scales as
\[
M = \Theta\!\left( \frac{N^2 \log N}{B} \right).
\]

For such a modulus, the pair sum generation cost scales as 
\[
\Theta(N^2 \log N \cdot M) = \Theta\!\left(\frac{N^4 \log^2 N}{B}\right).
\]
For sorting, each sum requires $\Theta(\log N)$ bits of storage, meaning each subproblem accommodates $\Theta(B / \log N)$ sums.
Sorting these sums within a subproblem scales as $\Theta\bigl( \frac{B}{\log N}\log(\frac{B}{\log N}) \log N \bigr)=\Theta\bigl(B\log(\frac{B}{\log N})\bigr)$.
Since the search space is partitioned across $\Theta(M^2)$ subproblems, the total sorting cost across all subproblems scales as $\Theta\bigl(M^2 B\log(\frac{B}{\log N})\bigr)$.
Using $M = \Theta(N^2 \log N / B)$, this becomes
\[
\Theta\!\left(
\frac{N^4 \log^2 N}{B}
\log\!\left(\frac{B}{\log N}\right)
\right).
\]
Assuming $B \gg \log N$, the sorting phase dominates the combined preprocessing complexity.

Evaluation cost remains $\Theta(N^3 \log N)$ under partitioning.
As a result, the preprocessing-to-evaluation ratio, which without partitioning is
\[
\frac{
\Theta(N^2 \log^2 N)
}{
\Theta(N^3 \log N)
}
=
\Theta\!\left(\frac{\log N}{N}\right)
\longrightarrow 0
\quad \text{as } N \to \infty,
\]
becomes under fixed-memory partitioning
\[
\frac{
\Theta\!\left(
\frac{N^4 \log^2 N}{B}
\log\!\left(\frac{B}{\log N}\right)
\right)
}{
\Theta(N^3 \log N)
}
=
\Theta\!\left(
\frac{
N \log N
}{
B
}
\log\!\left(\frac{B}{\log N}\right)
\right)
\longrightarrow \infty
\quad \text{as } N \to \infty.
\]

Thus, preprocessing eventually dominates runtime as $e_{\max}$ grows under a fixed memory budget, with the crossover delayed by increasing $B$.
Empirically, preprocessing remains a small fraction of runtime across our search ranges, as shown in Section~\ref{sec:results}.

Preprocessing dominance intensifies as the number of terms increases.
A search for solutions to the seven-term case $\sum_{i=1}^{6} a_i^n = a_7^n$ requires storing three-term sums $x_1^n + x_2^n + x_3^n$, giving an unpartitioned search space of size $\Theta(N^3 \log N)$ bits.
Under the same memory budget $B$, the required modulus scales as $M = \Theta(N^3 \log N / B)$.
Accounting for both pair sum generation and sorting, the total preprocessing cost scales as $\Theta\bigl( \frac{N^6 \log^2 N}{B}\log(\frac{B}{\log N}) \bigr)$, while evaluation scales as $\Theta(N^4 \log N)$.
Consequently, the preprocessing-to-evaluation ratio becomes $\Theta\bigl( \frac{N^2 \log N}{B}\log(\frac{B}{\log N}) \bigr)$, diverging more rapidly than in the five-term case.

More generally, for a $k$-term equation where the left-hand side is decomposed into two sums of $t = (k-1)/2$ terms, the preprocessing-to-evaluation ratio scales as $\Theta\bigl( \frac{N^{t-1} \log N}{B} \log\bigl(\frac{B}{\log N}\bigr) \bigr)$.

\subsection{Reducing the Search Space}
\label{subsec:reducing_search_space}

We reduce the search space using bounding constraints and modular residue filtering.

When generating sums $x_1^5 + x_2^5$, we discard pairs with $x_1^5 + x_2^5 \ge e_{\max}^5$, retaining only candidates that can contribute to solutions in $\mathbb{Z}_{\ge 0}$.

Stronger reductions arise from congruence constraints.
For many moduli $m$, the image of $x^5 \bmod m$ is sparse.
For example, modulo $11$,
\[
x^5 \equiv \{0,\pm 1\} \pmod{11}
\implies
x_1^5 + x_2^5 \equiv \{0, \pm 1, \pm 2\} \pmod{11}.
\]
Hence admissible pairs depend only on $e^5 \bmod 11$; e.g., if $e^5 \equiv 1 \pmod{11}$, residues summing to $-2$ are excluded, and if $e^5 \equiv -1 \pmod{11}$, residues summing to $2$ are excluded.
This yields three residue classes $e^5 \equiv 0, \pm 1 \pmod{11}$, with the non-zero cases allowing iteration over a smaller multiset $S$.
The three resulting multisets are not disjoint; the same sum $x_1^5 + x_2^5$ may appear in many multisets and is generated and sorted independently for each.

The same principle applies for $m \in \{25,31,41,61\}$; Table~\ref{tab:residue_classes} summarizes residue class counts and resulting search-space reductions for each $m$.

\begin{table}[htbp]
\centering
\caption{Iteration reduction from modular residue filtering.}
\label{tab:residue_classes}
\renewcommand{\arraystretch}{1.1}
\begin{tabular}{c c c}
\toprule
Modulus $m$ & Residue Classes & Reduction (\%) \\
\midrule
11 & 3  & 19 \\
25 & 5  & 22 \\
31 & 7  & 18 \\
41 & 9  & 17 \\
61 & 13 & 16 \\
\bottomrule
\end{tabular}
\end{table}

Combining moduli multiplies the number of partitions; using the full set would yield $3 \times 5 \times 7 \times 9 \times 13 = 12{,}285$ overlapping classes.
We restrict our implementation to $m=11$ and $25$, balancing search space reduction with repeated generation and sorting cost.
Additionally, when the partitioning modulus $M$ is divisible by $5$, the modulo-$25$ filter can be applied without a corresponding five-fold increase in preprocessing: since $x^5 \bmod 25$ depends only on
$x \bmod 5$, the partition already determines each summand's fifth-power residue modulo $25$.

\subsection{Extension to \(\mathbb{Z}\)}

To extend the algorithm to all integers, we apply three modifications.

\begin{enumerate}[label=(\roman*)]
\item Expand the enumeration of pairs in $S$ to
\[
S = \big( x_1^5 + x_2^5 \big)_{-e_{\max} \le x_1 \le x_2 \le e_{\max}},
\]
omitting the pruning condition $x_1^5 + x_2^5 < e_{\max}^5$ used for $\mathbb{Z}_{\geq 0}$.

\item Adjust the two-pointer initialization to
\[
s_1 = \min\{ s \in S : s \ge -e^5 \}, \quad s_2 = \max\{ s \in S : s < 2e^5 \},
\]
capturing the full range of admissible sums in $\mathbb{Z}$.

\item Discard pairs with a zero term to avoid trivial families of solutions (e.g., $0^5 + u^5 + (-u)^5 + v^5 = v^5$ for all integers $|u|, |v| \le e_{\max}$).
These trivial families generate $O(N^2)$ solutions that are expensive to detect in our implementation, as detailed in Section~\ref{sec:primitivity}.
    
To cover cases with a zero term, we implement a stream-merge search as in \cite{Bernstein2001} to exhaustively search the families $a^5 + b^5 + c^5 = d^5$ and $a^5 + b^5 = c^5 + d^5$ up to $e_{\max}$.
This auxiliary search operates at $O(N^2 \log^2 N)$ time complexity rather than $O(N^3 \log N)$ for the five-term pipeline.

We additionally partition the stream-merge search as in Section~\ref{subsec:partitioning}, using a modulus $M$ for which $x \mapsto x^5 \bmod M$ is injective, filtering the left-hand side and right-hand side to the same residue class.
This yields $M$ independent subproblems, here motivated by multithreading rather than memory constraints.

\end{enumerate}

These extensions from \(\mathbb{Z}_{\geq 0}\) to $\mathbb{Z}$ keep asymptotic complexity unchanged, while increasing the iteration count by a factor of approximately $4.8$.

\section{Implementation}

The algorithm was implemented in C++ using \texttt{g++}.
All values involved in fifth-power arithmetic used 128-bit signed integers (\verb|__int128|), sufficient for our search ranges.
The implementation relied on the compiler for vectorization, as the data-dependent two-pointer traversal precludes explicit SIMD.

\subsection{Instance Selection}

Due to memory scaling considerations described in Section~\ref{subsec:fixed_memory}, the search predominantly used hosts with 384 GiB of memory and 96 vCPUs.
While large, these are widely-available instance types on commercial cloud platforms.

\subsection{Difference Encoding}

Before each iteration over $S$, the array was transformed in-place into a difference-encoded form,
\[
d_i = S_i - S_{i-1}, \qquad d_0 = 0.
\]
This allowed incremental updates of $D = e^5 - (s_1 + s_2)$, where pointer advances modify $D$ by $\pm d_i$, eliminating one arithmetic operation per iteration and reducing runtime by approximately 25\%.

\subsection{Sorting}

The multiset $S$ was sorted using the in-place parallel IPS$^4$o algorithm \cite{Axtmann2022}, which outperformed \texttt{g++}'s parallel \texttt{std::sort}. In-place sorting was essential, as peak memory required by out-of-place methods would remain unused during the evaluation phase.

\subsection{Multithreading and Load Balancing}

All compute-intensive stages (generation of $S$, sorting, difference encoding, and iteration over $e$) were parallelized using OpenMP.

During iteration, parallelizing by assigning a single target $e$ to each thread led to load imbalance as the number of remaining targets fell below the thread count.
We addressed this by further partitioning the $s_2$ search range for each target $e$, improving thread utilization.

\subsection{Verification of Primitivity}
\label{sec:primitivity}

A practical challenge arises in verifying primitivity: when $D = 0$, the values $\{a,b,c,d\}$ are expensive to recover, and multiples of known solutions generate large numbers of matches, predominantly from the $e = 144$ solution given its small size.
We therefore apply a hierarchical screening procedure that first checks whether $e$ is a multiple of a known solution target (or, in $\mathbb{Z}$, whether $e$ matches any term from a known solution, since rearrangement can shift the target term), then verifies consistency with the known solution's pair sums.
Candidates passing both checks are fully reconstructed and reported; if imprimitive, dividing by $\gcd(a,b,c,d,e)$ yields a new primitive solution.
Reporting such imprimitive solutions facilitates partial-range searching, allowing a primitive solution missed in one range to be recovered later, as a multiple detected in a subsequent range.

\section{Results}
\label{sec:results}

The search combined exhaustive and partial coverage over the ranges in Table~\ref{tab:search_ranges}, where partial coverage denotes that a subset of partitions has been run.
This approach recovered all known solutions and identified the new solution reported here.

\begin{table}[htbp]
\centering
\caption{Search ranges for primitive solutions.}
\label{tab:search_ranges}
\renewcommand{\arraystretch}{1.2}
\begin{tabular}{c c c}
\toprule
Domain & Exhaustive Coverage & Partial Coverage \\
\midrule
$\mathbb{Z}_{\ge 0}$ & $e \le 2.76 \times 10^6$ & $2.76 \times 10^6 < e \le 4.00 \times 10^6$ \\
$\mathbb{Z}$          & $e \le 1.45 \times 10^6$ & $1.45 \times 10^6 < e \le 2.00 \times 10^6$ \\
\bottomrule
\end{tabular}
\end{table}

Exhaustive coverage was verified by ensuring the processed target count matched the expected total over each search range.
Since each target residue class $e \pmod M$ corresponds to $(M+1)/2$ unique residue pairs $\{r_1, r_2\} \pmod M$, the expected number of evaluated targets across all partitions is given by
\[
(e_{\max} - e_{\min} + 1)\,\frac{M+1}{2}.
\]
All completed ranges satisfied this identity.

All computations were performed using error-correcting memory.
However, at our scale, both soft errors (e.g., bit flips) and hardware faults cannot be fully excluded.
To provide robustness against such errors, for a solution whose four summands are distinct, our algorithm detects it exactly three times in $\mathbb{Z}_{\geq 0}$ and at least three times in $\mathbb{Z}$.
For example, the solution with $e = 144$ is detected as:
\begin{align*}
(s_1, s_2) &= (27^5 + 84^5,\; 110^5 + 133^5),\\
           &= (27^5 + 110^5,\; 84^5 + 133^5),\\
           &= (27^5 + 133^5,\; 84^5 + 110^5).
\end{align*}
The only candidate solutions arising without multiplicity would have \(a=b=c=d\), which yields \(4a^5=e^5\), impossible since $4$ is not a perfect fifth power.
All known solutions were recovered with the expected counts. 

Regression testing against a baseline implementation verified that optimizations preserved exact iteration counts.

The complete search encompassed a total of $5.4 \times 10^{18}$ iterations (evaluations of the target condition $D=0$), requiring approximately 10.5 million vCPU-hours over a nine-month period.

Table~\ref{tab:preprocessing_ratio} reports the preprocessing-to-evaluation ratio across benchmark runs in $\mathbb{Z}$, with $e_{\min} = 0.75e_{\max}$, 384 GiB of memory, and $M$ chosen so that $S$ occupies 85\% of available memory.
These ratios reflect the partitioning scheme of Section~\ref{subsec:partitioning} in isolation; modular residue filtering (Section~\ref{subsec:reducing_search_space}) is excluded, as its repeated generation and sorting overhead becomes increasingly disadvantageous relative to evaluation as search ranges grow, consistent with the preprocessing dominance analyzed in Section~\ref{subsec:fixed_memory}.
Even where preprocessing begins to dominate, the total runtime remains favorable relative to lower-memory alternatives that bypass preprocessing entirely, namely the priority-queue search, as our benchmarking (Section~\ref{subsec:alt_algorithms}) showed the two-pointer search to be more than an order of magnitude faster.

\begin{table}[htbp]
\centering
\caption{Preprocessing-to-evaluation ratio in $\mathbb{Z}$, with $e_{\min} = 0.75e_{\max}$; excludes modular residue filtering (Section~\ref{subsec:reducing_search_space}).}
\label{tab:preprocessing_ratio}
\renewcommand{\arraystretch}{1.1}
\begin{tabular}{c c c}
\toprule
$e_{\max}$ & $M$ & Ratio \\
\midrule
$2 \times 10^6$ & 749 & 0.065 \\
$4 \times 10^6$ & 2999 & 0.127 \\
$8 \times 10^6$ & 11998 & 0.248 \\
$1.6 \times 10^7$ & 47999 & 0.500 \\
$3.2\times 10^7$ & 191999 & 1.018 \\
\bottomrule
\end{tabular}
\end{table}

\section{Distribution of Terms and Comparison with \(n=4\)}
\label{sec:distribution}

The discovery of a fourth primitive solution to \eqref{eq:main} enables a preliminary comparison between counterexamples to Euler's sum-of-powers conjecture for \(n=4\) and \(n=5\).
This comparison was motivated by an observation of Noam Elkies that the two largest terms in \eqref{eq:new_solution} are nearly equal, reminiscent of the near-equality of two parameters in Booker and Sutherland's solution for $3$ as a sum of three cubes \cite{BookerSutherland2021}.

For a primitive solution to \eqref{eq:sum_of_powers}, define the dominance ratio
\[
\rho =
\left(
\frac{
\max_{1 \le i < k} |a_i|
}{
|a_k|
}
\right)^n.
\]
Values of $\rho$ near $1$ indicate that the largest summand is close in magnitude to the target.

Figure~\ref{fig:ratios_plot} compares known primitive solutions for \(n=4\) and \(n=5\).
The \(n=4\) data consist of the primitive solutions found in the exhaustive search of Gerbicz, Durman, Radaev, and Zubkov \cite{Gerbicz2009} over the reported range.
The \(n=5\) data consist of the four currently known primitive solutions, all recovered by our search.
Our partitioning scheme is independent of relative term magnitudes, so partial coverage does not preferentially favor large or small values of $\rho$.

\begin{figure}[htbp]
\centering
\begin{tikzpicture}
\begin{axis}[
    width=0.75\textwidth,
    height=6.8cm,
    xlabel={Target value \(|a_k|\)},
    ylabel={Dominance ratio \(\rho\)},
    xmode=log,
    xmin=100,
    xmax=2.5e9,
    ymin=0.45,
    ymax=1.02,
    axis lines=box,
    tick label style={font=\small},
    label style={font=\small},
    legend style={
        font=\small,
        draw=none,
        fill=none
    },
    legend pos=outer north east
]

\addplot[
    only marks,
    mark=o,
    mark size=1.7pt
] coordinates {
    (422481, 0.9271)
    (2813001, 0.9370)
    (8707481, 0.8384)
    (12197457, 0.7338)
    (16003017, 0.6154)
    (16430513, 0.9641)
    (20615673, 0.6911)
    (44310257, 0.7391)
    (68711097, 0.8478)
    (117112081, 0.6752)
    (145087793, 0.5008)
    (156646737, 0.7677)
    (589845921, 0.9522)
    (638523249, 0.9517)
    (873822121, 0.6008)
    (1259768473, 0.7357)
    (1679142729, 0.9798)
    (1787882337, 0.7485)
    (1871713857, 0.5254)
};
\addlegendentry{\(n=4\)}

\addplot[
    only marks,
    mark=*,
    mark size=2pt
] coordinates {
    (144, 0.6721)
    (14132, 0.9776)
    (85359, 0.9955)
    (1956878, 0.9983)
};
\addplot[dashed, thin, domain=100:2.5e9] {1};
\addlegendentry{\(n=5\)}

\end{axis}
\end{tikzpicture}

\caption{Dominance ratios \(\rho\) of known primitive solutions for \(n=4\) and \(n=5\).}
\label{fig:ratios_plot}
\end{figure}

As shown in Figure~\ref{fig:ratios_plot}, the known primitive solutions for \(n=4\) exhibit substantial variation in \(\rho\), with median \(\approx 0.7485\).

In contrast, while the smallest known \(n=5\) solution has \(\rho \approx 0.6721\), the remaining three all have \(\rho > 0.97\), with the new solution \eqref{eq:new_solution} reaching \(\rho \approx 0.9983\).

We assess whether this clustering is anomalous using a geometric baseline model in which the summands are drawn independently and uniformly at random.
Specifically, we draw $k-1$ summands uniformly from $[0,1]^{k-1}$ (for $\mathbb{Z}_{\geq 0}$ solutions) or $[-1,1]^{k-1}$ (for $\mathbb{Z}$ solutions), then compute the implied $k$th term as the real $n$th root of their signed power sum. We run $10^7$ independent trials per configuration.
For these draws, we compute $\rho$ as the $n$th power of the ratio of the two largest absolute values among all $k$ terms, rather than assuming the last term is the largest. This matches the definition above for known solutions, since they are conventionally arranged so that $|a_k|$ is the largest term.

The nineteen known $n=4$ solutions track their baseline distribution without anomalous clustering, with an empirical median at the 46th percentile.
Among the $n=5$ solutions, the 1966 and 1996 solutions fall at the 44th and 87th percentiles of their respective baselines, while the 2004 and new solutions each sit at the 97th percentile.
We note that the 1966 solution's integrality caps its $\rho \le (143/144)^5 \approx 0.9657$, a constraint absent from the continuous baseline, which may partly explain its comparatively low percentile in the baseline model.

With only four data points, the concentration of three solutions in the upper tail may reflect an underlying property of fifth-power solutions, but is far from conclusive.
This uncertainty highlights the value of continued computational search for further solutions: as the corpus of solutions grows, structural properties can be studied with greater confidence.

\section{Potential Further Optimizations}

One additional refinement, not implemented in our algorithm, would multiplicatively combine all moduli in Table~\ref{tab:residue_classes} as
\[
\tilde{M} = 11 \cdot 25 \cdot 31 \cdot 41 \cdot 61.
\]

Instead of separately generating $S$ for each residue-class combination across the moduli in Table~\ref{tab:residue_classes}, all sums $x_1^5 + x_2^5$ could be computed once and partitioned by their residues modulo $\tilde{M}$, with each partition sorted independently.

For each target value $e$, only partitions satisfying
\[
s_1 + s_2 \equiv e^5 \pmod{\tilde{M}}
\]
would be processed during iteration.

This would allow full residue class filtering without repeated generation and sorting of $S$. In principle, this could yield an additional factor-of-two speedup.

\section*{Code Availability}
Code to reproduce our results is available at:
\begin{center}
\url{https://github.com/JeffreyMBraun/Euler514}
\end{center}

\section*{Acknowledgments}

The author gratefully acknowledges the support and encouragement of his spouse, Randie Kim, and his family throughout this work.

The author thanks Noam Elkies for the observation that motivated Section~\ref{sec:distribution}, and Drew Sutherland for communicating it and for helpful feedback on an earlier version of this manuscript.

\bibliographystyle{amsalpha}
\bibliography{references}

@article{Lander1966,
  author  = {Lander, L.~J. and Parkin, T.~R.},
  title   = {Counterexample to {E}uler's conjecture on sums of like powers},
  journal = {Bull. Amer. Math. Soc.},
  volume  = {72},
  year    = {1966},
  pages   = {1079}
}

@article{Lander1967,
  author  = {Lander, L.~J. and Parkin, T.~R.},
  title   = {A counterexample to {E}uler's sum of powers conjecture},
  journal = {Math. Comp.},
  volume  = {21},
  year    = {1967},
  pages   = {101--103}
}

@misc{Scher1996,
  author       = {Scher, B. and Seidl, E.},
  title        = {Massively Parallel Number Theory},
  howpublished = {Research announcement},
  year         = {1996},
  note         = {Archived at \url{https://web.archive.org/web/19970608041930/http://midway.ca.sandia.gov/~mecolv/euler/}}
}

@misc{Frye2004,
  author = {Frye, J.},
  title  = {A solution to $a^5 + b^5 + c^5 + d^5 = e^5$},
  year   = {2004},
  note   = {Unpublished; reported in I. Peterson, ``Euler's Sums of Powers,'' \textit{Science News}, September 24, 2004; see also J.-C. Meyrignac, pers. comm. to \textit{Wolfram MathWorld}, September 9, 2004.}
}

@article{Axtmann2022,
  author  = {Axtmann, M. and Witt, S. and Ferizovic, D. and Sanders, P.},
  title   = {Engineering In-place (Shared-memory) Sorting Algorithms},
  journal = {{ACM} Trans. Parallel Comput.},
  volume  = {9},
  number  = {1},
  year    = {2022},
  pages   = {1--62},
  doi     = {10.1145/3505181},
  note    = {Implementation available at \url{https://github.com/ips4o/ips4o}}
}

@article{Bernstein2001,
  author  = {Bernstein, D.~J.},
  title   = {Enumerating solutions to $p(a)+q(b)=r(c)+s(d)$},
  journal = {Math. Comp.},
  volume  = {70},
  number  = {233},
  year    = {2001},
  pages   = {389--394},
  doi     = {10.1090/S0025-5718-00-01219-9}
}

@article{Elkies1988,
  author  = {Elkies, N.~D.},
  title   = {On ${A}^4 + {B}^4 + {C}^4 = {D}^4$},
  journal = {Math. Comp.},
  volume  = {51},
  number  = {184},
  year    = {1988},
  pages   = {825--835},
  doi     = {10.1090/S0025-5718-1988-0935042-8}
}

@misc{Gerbicz2009,
  author       = {Gerbicz, R. and Durman, L. and Radaev, Y. and Zubkov, A.},
  title        = {{E}uler {Q}uartic {C}onjecture: ${A}^4 = {B}^4 + {C}^4 + {D}^4$},
  howpublished = {Webpage},
  year         = {2009},
  note         = {Archived at \url{https://web.archive.org/web/20230418025309/http://euler413.narod.ru/}}
}

@article{BookerSutherland2021,
author = {Booker, A.~R. and Sutherland, A.~V.},
title = {On a question of {M}ordell},
journal = {Proc. Natl. Acad. Sci. USA},
volume = {118},
number = {11},
year = {2021},
pages = {e2022377118},
doi = {10.1073/pnas.2022377118}
}

\end{document}